# NORMALIZED RANDOM MEASURES DRIVEN BY INCREASING ADDITIVE PROCESSES

By Luis E. Nieto-Barajas[1], Igor Prünster[2] and Stephen G. Walker[3]

*ITAM-México, Università degli Studi di Pavia and University of Kent*

This paper introduces and studies a new class of nonparametric prior distributions. Random probability distribution functions are constructed via normalization of random measures driven by increasing additive processes. In particular, we present results for the distribution of means under both prior and posterior conditions and, via the use of strategic latent variables, undertake a full Bayesian analysis. Our class of priors includes the well-known and widely used mixture of a Dirichlet process.

**1. Introduction.** This paper considers the problem of constructing a stochastic process, defined on the real line, which has sample paths behaving almost surely (a.s.) as a probability distribution function (d.f.). The law governing the process acts as a prior in Bayesian nonparametric problems. One popular idea is to take the random probability d.f. as a normalized increasing process $F(t) = Z(t)/\bar{Z}$, where $\bar{Z} = \lim_{t\to\infty} Z(t) < +\infty$ (a.s.). For example, the Dirichlet process [Ferguson (1973)] arises when $Z$ is a suitably reparameterized gamma process. We consider the case of *normalized random d.f.s driven by an increasing additive process* (IAP) $L$, that is, $Z(t) = \int k(t,x)\,dL(x)$ and provide regularity conditions on $k$ and $L$ to ensure $F$ is a random probability d.f. (a.s.).

This paper represents a natural development of the work of Regazzini, Lijoi and Prünster (2003). These authors introduce the class of normalized

Received August 2002; revised July 2003.
[1]Supported in part by a CONACYT (Mexico) and an ORS (UK) grant.
[2]Supported in part by MIUR Young Researchers Program and by a CNR grant, while visiting the Department of Mathematical Sciences, University of Bath, UK.
[3]Supported by an EPSRC Advanced Research Fellowship.

random d.f.s with independent increments (RMI), a particular case of IAP driven random d.f.s, and consider the problem of determining the exact distribution of means of a normalized RMI. The study of means of random probability d.f.s has become, after the pioneering work of Cifarelli and Regazzini (1979, 1990), a very active area of research, touching on both analytical and simulation based approaches. We mention, among others, Muliere and Tardella (1998), Guglielmi and Tweedie (2001), Regazzini, Guglielmi and Di Nunno (2002) and Lijoi and Regazzini (2004) for results in the Dirichlet case and Epifani, Lijoi and Prünster (2003), Hjort (2003) and Regazzini, Lijoi and Prünster (2003) for results beyond the Dirichlet case.

In this paper we also consider both analytical and simulation based approaches to the study of means, providing a comprehensive treatment of the subject. By extending the methodology proposed in Regazzini, Lijoi and Prünster (2003) to our more general case, we determine the exact law of any mean of a normalized IAP driven random d.f. This approach exploits Gurland's inversion formula and gives expressions for the posterior distributions in terms of the Liouville–Weyl fractional integral.

An important class of normalized IAP driven random d.f.s is obtained if $\lim_{t\to\infty} k(t,x)$ is a constant for all $x$. Then we obtain the class of mixtures of normalized RMI, that is, $F(t)=\int k(t,x)\,dG(x)$. Moreover, if $G$ is a Dirichlet process we have a mixture of a Dirichlet process, first introduced by Lo (1984). This family was the focus of much attention during the 1990s as a consequence of the introduction of simulation based inference, first considered by Escobar (1988) in his Ph.D. dissertation, and later developed by Escobar and West (1995) and MacEachern and Müller (1998), among others. The model is comprehensively reviewed in the book edited by Dey, Müller and Sinha (1998). By exploiting the above mentioned general results, we are able to give exact prior and posterior distributions for any mean of a mixture of a Dirichlet prior and, furthermore, provide a new simulation algorithm. Finally, we illustrate our results both theoretically and numerically by applying them to what we call a Dirichlet driven random probability d.f.

Before proceeding we introduce the fundamental concepts and tools for the paper. Let $L := \{L(y): y \ge 0\}$ be any IAP defined on $(\Omega, \mathcal{F}, P)$. In general, an IAP can be expressed as

$$L(y) = \sum_{j\,:\,\tau_j \le y} L\{\tau_j\} + L^c(y), \tag{1}$$

where $M = \{\tau_1, \tau_2, \ldots\}$ is the set of fixed points of discontinuity and $L^c$ is the part of the process without fixed points of discontinuity. Hence $L$ is characterized by the density functions of the jumps $\{L\{\tau_1\}, L\{\tau_2\}, \ldots\}$, indicated by $\{f_{\tau_1}, f_{\tau_2}, \ldots\}$, and the family of Lévy measures, $\{\nu : \nu(y, dv), y \ge 0\}$, related to $L^c$ through the celebrated Lévy–Khintchine representation.

For an exhaustive account of the theory of IAPs see, for example, Sato (1999).

Consider any nondegenerate measure $\alpha$ on $\mathscr{B}(\mathbb{R})$ such that $\alpha(\mathbb{R}) = a \in (0, +\infty)$ and denote by $A$ the corresponding d.f. The time change $y = A(x)$ yields an a.s. $[P]$ finite IAP $L_A = \{L_A(x) : x \in \mathbb{R}\}$ uniquely determined by the family of Lévy measures $\{\nu_\alpha : \nu_\alpha(x, dv) = \nu(G_\alpha^{-1}(x), dv), x \in \mathbb{R}\}$, where $G_\alpha(y) := \inf\{x : A(x) \geq y\}$ for $y \in (0, \alpha(\mathbb{R}))$. Its Laplace transform is thus given by

$$E[e^{-\lambda L_A(x)}] = \exp\left[-\int_0^\infty \{1 - e^{-\lambda v}\} \nu_\alpha(x, dv)\right] \qquad \text{for any } \lambda \geq 0.$$

In the following it will be convenient to characterize $L_A$ in terms of its Poisson intensity measure, indicated by $\tilde{\nu}_\alpha$, instead of its corresponding family of Lévy measures. Recall that $\tilde{\nu}_\alpha((-\infty, x] \times C) = \nu_\alpha(x, C)$ for every $x \in \mathbb{R}$ and $C \in \mathscr{B}((0, +\infty))$.

Consider now the stochastic process given by a convolution of $L_A$ with $k : \mathbb{R} \times \mathbb{R} \to \mathbb{R}^+$,

$$Z := \left\{ Z(t) = \int_{\mathbb{R}} k(t, x)\, dL_A(x) : t \in \mathbb{R} \right\}.$$

Suppose $k$ and $\tilde{\nu}_\alpha$ satisfy simultaneously the following conditions:

(I) $t \mapsto k(t, x)$ is nondecreasing and right continuous with $\lim_{t \to -\infty} k(t, x) = 0$ for every $x \in \mathbb{R}$;
(II) $\int_{\mathbb{R} \times (0, +\infty)} [1 - \exp\{-\lambda v \bar{k}(x)\}] \tilde{\nu}_\alpha(dx\, dv) < +\infty$ for every $\lambda > 0$, where $\bar{k}(x) := \lim_{t \to +\infty} k(t, x)$;
(III) $\tilde{\nu}_\alpha(\mathbb{R} \times (0, +\infty)) = +\infty$.

Then $Z$ is a random d.f. a.s. $[P]$ and $F = \{F(t) = Z(t)/\bar{Z} : t \in \mathbb{R}\}$ is a random probability d.f. a.s. $[P]$, having set $\bar{Z} := \lim_{t \to +\infty} Z(t)$. For details about the determination of conditions (I)–(III), refer to the Appendix.

In this context, according to Barndorff-Nielsen and Shephard (2001), $L_A$ can be seen as a *background driving IAP*. Hence, $Z$ and $F$ will be called an *IAP driven random d.f.* and a *normalized IAP driven random d.f.*, respectively. There are now a number of works based on Lévy driven processes; we mention Wolpert, Ickstadt and Hansen (2003) and Brockwell (2001) who introduce a Lévy driven CARMA model. By choosing $k(t, x) = \mathbb{I}_{(-\infty, t]}(x)$ a normalized IAP driven random d.f. reduces to a normalized RMI, whose trajectories are discrete a.s. $[P]$. This property of normalized RMI may be undesirable in many situations. It is easily seen that a normalized IAP driven random d.f. has absolutely continuous sample paths with respect to the Lebesgue measure on $\mathbb{R}$ a.s. $[P]$ if and only if $t \mapsto k(t, x)$ is absolutely continuous for every $x \in \mathbb{R}$. If this is the case, the corresponding normalized IAP driven random density function is given by

$$f(t) = \frac{\int_{\mathbb{R}} k'(t,x)\, dL_A(x)}{\bar{Z}}, \qquad (t \in \mathbb{R}) \text{ a.s. } [P],$$

where $k'(t,x) := \frac{\partial}{\partial t} k(t,x)$. In the following we will always assume $F$ to admit a density.

In Section 2 we derive the exact distributions of means of normalized IAP driven random measures under prior and posterior conditions and derive distributional results for means of normalized gamma driven random d.f.s and, in particular, for the mixture of a Dirichlet process. In Section 3 a sampling strategy for drawing samples from the posterior distribution of $F$ is presented, and we provide a numerical illustration. All proofs are deferred to the Appendix.

**2. Distribution of means of normalized IAP driven d.f.s.** In this section we are concerned with the problem of determining the prior and posterior distribution of means of normalized IAP driven random d.f.s, extending the results of Regazzini, Lijoi and Prünster (2003) (RLP).

2.1. *Existence and distribution of means.* First of all we need to establish the existence of $\int_{\mathbb{R}} g(t)\, dF(t)$, or equivalently of $\int_{\mathbb{R}} g(t)\, dZ(t)$. Suppose $\int_{\mathbb{R}} |g(t)|\, dZ(t) < +\infty$ a.s. $[P]$. By application of Fubini's theorem,

$$\begin{aligned} \int_{\mathbb{R}} g(t)\, dZ(t) &= \int_{\mathbb{R}} g(t) \int_{\mathbb{R}} k'(t,x)\, dL_A(x)\, dt \\ &= \int_{\mathbb{R}} h(x)\, dL_A(x) \qquad \text{a.s. } [P], \end{aligned} \tag{2}$$

where $h(x) = \int_{\mathbb{R}} g(t) k'(t,x)\, dt$. Hence, a linear functional of an IAP driven random d.f. can be expressed as another linear functional of an IAP, which actually reduces our problem to the one considered by RLP. In terms of existence, the previous relation guarantees that $\int_{\mathbb{R}} |g(t)|\, dF(t) = \int_{\mathbb{R}} \tilde{h}(x)\, dL_A(x)$ a.s. $[P]$, having set $\tilde{h}(x) = \int_{\mathbb{R}} |g(t)| k'(t,x)\, dt$. Thus, by a slight modification of Proposition 1 in RLP, we have the required necessary and sufficient condition.

PROPOSITION 1. *Let $F$ be any normalized IAP driven random d.f. and let $g$ be any measurable function $g : \mathbb{R} \to \mathbb{R}$. Set $\tilde{h}(x) := \int_{\mathbb{R}} |g(t)| k'(t,x)\, dt$. Then $\int_{\mathbb{R}} |g(t)|\, dF(t) < +\infty$ a.s. $[P]$ if and only if $\int_{\mathbb{R}\times(0,+\infty)} [1 - \exp(-\lambda v \tilde{h}(x))] \times \tilde{\nu}_\alpha(dx\, dv) < +\infty$ holds for every $\lambda > 0$.*

We now proceed to determine the probability distribution of $\int_{\mathbb{R}} g(t)\,dF(t)$. Assuming the conditions of Proposition 1 hold, we observe that for any $\sigma \in \mathbb{R}$,

$$P\left\{\int_{\mathbb{R}} g(t)\,dF(t) \le \sigma\right\} = P\left\{\int_{\mathbb{R}} \{h(x) - \sigma\bar{k}(x)\}\,dL_A(x) \le 0\right\},$$

where $h(x) := \int_{\mathbb{R}} g(t)k'(t,x)\,dt$ and $\bar{k}(x) := \lim_{t\to+\infty} k(t,x)$. Hence, we are able to extend Proposition 2 in RLP, which is based on the inversion formula given in Gurland (1948), to our more general case, with obvious modifications.

Proposition 2. *Let $F$ be a normalized IAP driven random d.f., let $\mathbb{F}$ be the probability d.f. of $\int_{\mathbb{R}} g(t)\,dF(t)$ and set $h(x) = \int_{\mathbb{R}} g(t)k'(t,x)\,dt$. For every $\sigma \in \mathbb{R}$, we have*

$$\begin{aligned}\frac{1}{2}&[\mathbb{F}(\sigma) + \mathbb{F}(\sigma - 0)]\\ &= \frac{1}{2} - \frac{1}{\pi}\lim_{T\uparrow+\infty}\int_0^T \frac{1}{s}\exp\left\{\int_{\mathbb{R}\times(0,+\infty)}[\cos(sv(h(x) - \sigma\bar{k}(x))) - 1]\tilde{\nu}_\alpha(dx\,dv)\right\}\\ &\qquad\times\sin\left(\int_{\mathbb{R}\times(0,+\infty)}\sin\{sv(h(x) - \sigma\bar{k}(x))\}\tilde{\nu}_\alpha(dx\,dv)\right)ds.\end{aligned}$$

2.2. *Posterior distribution of means.* Here we aim at providing expressions for the posterior distribution of means of normalized IAP driven random d.f.s. This is done by introducing an appropriate sequence of nested partitions and by discretizing $F$ through the discretization of both $k$ and the space of observations. This construction guarantees the discretized posterior distribution of the mean to determine uniquely the limiting one, by a.s. convergence in distribution. Hence, we give an explicit expression for the posterior density of the discretized mean, which can be used as an approximation of the limiting one. In certain cases, once the Lévy measure is specified, it is also possible to derive an explicit representation of the limiting distribution.

Assume that $(\Omega, \mathcal{F}, P)$ also supports a sequence $T = (T_n)_{n\ge1}$ of exchangeable random variables. The first step consists in discretizing $F$. To this end, let us introduce a sequence of partitions $(\mathcal{P}_m)_{m\ge1}$ of $\mathbb{R}$, where $\mathcal{P}_m := \{A_{m,i} : i = 0, \dots, k_m + 1\}$, which satisfy the following properties:

(a) $\mathcal{P}_{m+1}$ is a refinement of $\mathcal{P}_m$.
(b) $\mathscr{B}(\mathbb{R})$ is generated by $\bigcup_{m\ge1}\sigma(\mathcal{P}_m)$, where $\sigma(\mathcal{P}_m)$ denotes the $\sigma$-algebra generated by $\mathcal{P}_m$.
(c) $\varepsilon_m := 2\max_{1\le i\le k_m}\operatorname{diam}(A_{m,i}) \downarrow 0$ (as $m \to +\infty$).

(d) $A_{m,0} = (-\infty, -R_m)$, $A_{m,i} = [t_{m,i}, t_{m,i+1})$ for $i = 1, \ldots, k_m - 1$, $A_{m,k_m} = [t_{m,k_m}, t_{m,k_m+1}]$, and $A_{m,k_m+1} = (R_m, +\infty)$, with $t_{m,1} = -R_m$, $t_{m,k_m+1} = R_m$ and $R_m > 0$ for any $m \ge 1$.

Now we have to select points $a_{m,i}$ in $A_{m,i}$ for $i = 1, \ldots, k_m$ and put $a_{m,0} = -R_m$ and $a_{m,k_m+1} = R_m$. Whenever the $r$th element, $T_r$, in the sample lies in $A_{m,i}$, it is as if we had observed $a_{m,i}$. The discretized random d.f. is defined as

$$F_m(t) := \sum_{\{j : a_{m,j} \le t\}} \frac{\int_{\mathbb{R}} [k(t_{m,j+1}, x) - k(t_{m,j}, x)]\, dL_A(x)}{\bar{Z}} \qquad \text{for every } t \in \mathbb{R} \tag{3}$$

with the conventions $k(t_{m,0}, x) = 0$ and $k(t_{m,k_m+2}, x) = \lim_{t \to +\infty} k(t, x) = \bar{k}(x)$.

Set $T^n = (T_1, \ldots, T_n)$ and denote by $\mu^n$ its distribution. The prior probability d.f. of a mean $\int g(t)\, d\widetilde{F}_m(t)$ is denoted by $\mathbb{F}_m(\cdot; g)$ and its posterior d.f. by $\mathbb{F}_{m,t^n}(\cdot; g)$. Having (3) at hand, it is easy to verify that the approximation result given in Proposition 4 of RLP holds true also in this quite different setup. Hence, for every $\sigma$ belonging to the set of continuity points of $\mathbb{F}_{t^n}$,

$$\lim_{m \to +\infty} \mathbb{F}^*_{m,t^n}(\sigma; g) = \mathbb{F}_{t^n}(\sigma; g) \qquad \text{a.s. } [\mu^n]. \tag{4}$$

Having derived $\mathbb{F}_m$ according to Proposition 2, one can see that Proposition 3 in RLP extends also to our more general case. Thus, supposing $(a, b)$ is an interval containing all the $g(a_{m,i})$'s and assuming interchangeability of the derivative with the integral, one has that the posterior density function of $\int g(t)\, dF_m(t)$, given $T^n = t^n$ with $n_{i_p} > 0$ terms set equal to $a_{m,i_p}$ ($p = 1, \ldots, q$) such that $\sum_p n_{i_p} = n$, is given by

$$\begin{aligned} \rho_{m,t^n} &= \frac{(-1)^n}{\mu^n(t^n)} \\ &\quad \times \frac{\partial^n}{\partial r_{i_1}^{n_{i_1}} \cdots \partial r_{i_q}^{n_{i_q}}} I_{a^+}^{n-1} \mathbb{F}_m(\sigma; r_0, \ldots, r_{k+1}) \Big|_{(r_0, \ldots, r_{k+1}) = (g(a_{m,0}), \ldots, g(a_{m,k+1}))}, \end{aligned} \tag{5}$$

where $I_{a^+}^n h(\sigma) = \int_a^\sigma \frac{(\sigma - u)^{n-1}}{(n-1)!} h(u)\, du$ is the Liouville–Weyl fractional integral, for $n \ge 1$, and $I_{a^+}^0$ represents the identity operator.

2.3. *Normalized gamma and mixtures of Dirichlet process.* Many Bayesian nonparametric priors are constructed via transformations of gamma processes. Hence, it seems natural to focus attention on normalized gamma driven random d.f.s. Here, a complete treatment of the distributional properties of means of normalized gamma driven random d.f.s is provided.

Before proceeding, let us briefly recall that a reparameterized gamma process, $\Gamma_A$, is characterized by a Poisson intensity measure of the type $\tilde{\nu}_\alpha(dx\,dv) := e^{-v}v^{-1}\,dv\,d\alpha(x)$, while the extended gamma process, $\Gamma_A^\beta$, introduced by Dykstra and Laud (1981), is characterized by a Poisson intensity measure of the type $\tilde{\nu}_\alpha(dx\,dv) := e^{-\beta(x)v}v^{-1}\,dv\,d\alpha(x)$, where $\beta$ is a nonnegative piecewise continuous function. These two IAPs are connected to the normalized gamma driven random d.f. and, in particular, to the mixture of Dirichlet process (MDP) through the following relations.

PROPOSITION 3. *Suppose $F$ is a normalized IAP driven random d.f. Then:*

(i) *If $L_A = \Gamma_A^\beta$, $F$ can be represented as a normalized gamma driven random d.f.s, that is,*

$$\frac{\int_{\mathbb{R}} k(t,x)\,d\Gamma_A^\beta(x)}{\int_{\mathbb{R}} \bar{k}(x)\,d\Gamma_A^\beta(x)} = \frac{\int_{\mathbb{R}} k(t,x)(\beta(x))^{-1}\,d\Gamma_A(x)}{\int_{\mathbb{R}} \bar{k}(x)(\beta(x))^{-1}\,d\Gamma_A(x)} \qquad a.s.\ [P].$$

(ii) *If $L_A = \Gamma_A$, $F$ can be represented as a mixture of a normalized extended gamma process, that is,*

$$\frac{\int_{\mathbb{R}} k(t,x)\,d\Gamma_A(x)}{\int_{\mathbb{R}} \bar{k}(x)\,d\Gamma_A(x)} = \int_{\mathbb{R}} \frac{k(t,x)}{\bar{k}(x)} \frac{d\Gamma_A^{1/\bar{k}}(x)}{\bar{\Gamma}_A^{1/\bar{k}}} \qquad a.s.\ [P],$$

*having set $\bar{\Gamma}_A^{1/\bar{k}} := \lim_{x\to+\infty} \Gamma_A^{1/\bar{k}}(x)$. Moreover, if $\bar{k}(x) = b^{-1}$, then $F$ is an MDP, $\int_{\mathbb{R}} bk(t,x)\,d\mathscr{D}_A(x)$, where $\mathscr{D}_A$ denotes the Dirichlet random d.f.*

(iii) *If $L_A = \Gamma_A$, any mean of $F$, provided it exists, may be represented as a mean of a normalized extended gamma process, that is,*

$$\int_{\mathbb{R}} g(t)\,dF(t) = \int_{\mathbb{R}} \bar{h}(x) \frac{d\Gamma_A^{1/\bar{k}}(x)}{\Gamma_A^{1/\bar{k}}} \qquad a.s.\ [P],$$

*with $\bar{h}(x) = (\bar{k}(x))^{-1}\int_{\mathbb{R}} g(t)k'(t,x)\,dt$. If, moreover, $\bar{k}(x) = b^{-1}$, then $\int_{\mathbb{R}} g(t)\,dF(t)$ becomes a mean of a Dirichlet process, $\int_{\mathbb{R}} bh(x)\,d\mathscr{D}_A(x)$, where $h(x) = \int_{\mathbb{R}} g(t)k'(t,x)\,dt$.*

Thus, we have that a normalized extended gamma driven random d.f. is equivalent to a normalized gamma driven random d.f. and that MDPs are a special case. Nonetheless, in studying means of normalized gamma driven random d.f.s we confine ourselves to MDPs, because to date nothing is known about exact distributions of their means. This is done without loss of generality, since the following results are easily extended to any normalized gamma driven random d.f.

With reference to existence of a mean of an MDP, $\int_{\mathbb{R}} g(t)\,dF(t)$, by Proposition 3, the condition reduces to the well-known

$$\int_{\mathbb{R}} \log(1+\lambda|h(x)|)\alpha(dx) < +\infty \qquad \text{for every } \lambda > 0, \tag{6}$$

with $h(x) = \int_{\mathbb{R}} g(t) k'(t,x)\,dt$ as previously. See Feigin and Tweedie (1989) and Cifarelli and Regazzini (1990, 1996). Consequently the d.f. of a mean $\int_{\mathbb{R}} g(t)\,dF(t)$ is given by

$$\begin{aligned}\mathbb{F}(\sigma) = \frac{1}{2} - \frac{1}{\pi}\int_0^{+\infty} \frac{1}{s} \exp\biggl\{-\int_{\mathbb{R}} \log\{1 + s^2(h(x)-\sigma)^2\}\alpha(dx)\biggr\} \\ \times \sin\biggl(\int_{\mathbb{R}} \arctan[s(h(x)-\sigma)]\alpha(dx)\biggr)\,ds.\end{aligned} \tag{7}$$

The fact that our mean is just another mean with respect to the Dirichlet process implies that its law is absolutely continuous with respect to the Lebesgue measure. See Regazzini, Guglielmi and Di Nunno (2002) for expressions of the corresponding density function.

We now move on in stating our main result, which provides intuitive insight into the mixing character of the posterior behavior of means of MDPs.

THEOREM 1. *Suppose $F$ is an MDP and its mean $\int_{\mathbb{R}} g(t)\,dF(t)$ exists, that is, $g$ satisfies* (6). *Then its posterior distribution, given $T^n = t^n$, is absolutely continuous* (*with respect to the Lebesgue measure on $\mathbb{R}$*) *and a posterior probability density function is given by*

$$\rho_{t^n}(\sigma) = \int_{\mathbb{R}^n} \rho_{u^n}(\sigma) G(du_1,\dots,du_n|t^n), \tag{8}$$

*where*

$$G(du_1,\dots,du_n|t^n) = \frac{\prod_{j=1}^n k'(t_j,u_j)\alpha^n(du_1,\dots,du_n)}{\int_{\mathbb{R}^n}\prod_{j=1}^n k'(t_j,u_j)\alpha^n(du_1,\dots,du_n)}$$

*represents the distribution of the latent variables $U^n$, given the observations $T^n = t^n$, with $\alpha^n$ defined as the n-fold product measure $\prod_{k=1}^n(\alpha + \sum_{i=1}^{k-1}\delta_{u_i})$, and $\rho_{u^n}$ denotes the posterior distribution of $\int_{\mathbb{R}} h(x)\,d\mathscr{D}_A(x)$, given $U^n = u^n$, with $h(x) = \int_{\mathbb{R}} g(t)k'(t,x)\,dt$, and given by*

$$\rho_{u^n}(\sigma) = \frac{a}{\pi}\int_0^{+\infty} \mathrm{Re}\biggl(\exp\biggl\{-\int_{\mathbb{R}} \log[1 + is(h(x)-\sigma)]\alpha^*(dx)\biggr\}\biggr)\,ds,$$

*having set $\alpha^* = \alpha + \sum_{i=1}^n \delta_{u_i}$.*

A deficiency of the previous intuitive result is represented by the dimension of the integration region in (8), which grows as the sample size grows. This can be overcome by an application of Lemma 2 in Lo (1984), which essentially allows one to account for coincidences within the latent observations. To this end let us introduce some notation. Denote by $\mathcal{P} := \{C_i : i = 1,\dots,N(\mathcal{P})\}$ a partition of $\{1,2,\dots,n\}$, where $N(\mathcal{P})$ indicates the number of cells and $C_i$ the $i$th cell in the partition. Moreover, let $c_i$ be the number of elements in $C_i$.

COROLLARY 1. *Suppose $\int_{\mathbb{R}} g(t)\,dF(t)$ is a mean of an MDP and $g$ satisfies (6). Then its posterior density function, given $T^n = t^n$, is given by*

$$\rho_{t^n}(\sigma) = \frac{\sum_{\mathcal{P}}(\prod_{i=1}^{N(\mathcal{P})}[(c_i-1)!\int_{\mathbb{R}} \rho_{u^{c_i}}(\sigma)\prod_{p\in C_i} k'(t_p,u)\alpha(du)])}{\sum_{\mathcal{P}}(\prod_{i=1}^{N(\mathcal{P})}[(c_i-1)!\int_{\mathbb{R}} \prod_{p\in C_i} k'(t_p,u)\alpha(du)])},$$

*where $\rho_{u^{c_i}}$ denotes the posterior density of $\int_{\mathbb{R}} h(x)\,d\mathscr{D}_A(x)$, given $c_i$ observations equal to $u$.*

It is worth pointing out that the burden involved in posterior densities, when dealing with more than a few observations, becomes overwhelming for currently available computational tools. The necessity of a simulation algorithm is evident.

**3. Posterior simulation.** In this section we provide a method to sample from the posterior distribution of $F$, and $f$, given a set of $n$ observations $T^n$. The algorithm depends on the strategic and novel introduction of latent variables. Let $S$ and $U$ be latent variables, and consider the joint distribution

$$p(t,s,u|L_A) = \exp(-u\bar{Z})k'(t,s)\,dL_A(s), \qquad u \ge 0, s \in \mathbb{R},$$

where, as previously, $L_A$ is a reparameterized IAP and $k'(t,s) = \frac{\partial}{\partial t}k(t,s)$.

Note that $L_A$ is a pure jump process and so the support of $s$ will be the location of the jumps of $L_A$, that is, $p(s|t,u,L_A) \propto L_A\{s\}k'(t,s)$ and $L_A\{s\} = L_A(s) - L_A(s-)$. Clearly $p(t|L_A) = f(t)$, as required.

Having established the general sampling strategy, let us consider, in particular, normalized gamma driven random d.f.s. For computational reasons, we allow $L_A$ to have fixed points of discontinuity. Recall the representation of such an IAP given in (1) together with the related notation. We work with a normalized extended gamma driven random d.f. which we know to be equivalent to a normalized gamma driven random d.f. by Proposition 3. Some other authors have obtained posterior distributions when working with additive processes in different contexts [see, e.g., Dykstra and Laud (1981), Hjort (1990), Walker and Muliere (1997) and Nieto-Barajas and Walker (2004)]. Let us start with a single observation $T_1$; then we obtain the following result, where $\mathcal{G}$ denotes a gamma distribution.

PROPOSITION 4. *Let $F(t) = Z(t)/Z(\Upsilon)$ be a normalized IAP driven random measure, where $\Upsilon$ is the maximum time up to where the process is observed and $T_1$ is a random sample from $F$. Denote by $M$ the set of prior fixed points of discontinuity of $L_A$ and by $\star$ an updated parameter/function.*

(i) *Given $T_1 = t_1$, $S_1 = s_1 \in M$ and $U_1 = u_1$, the posterior parameters are $M^\star = M$,*

$$f^\star_{\tau_j}(x) \propto \begin{cases} x e^{-u_1 k(\Upsilon,\tau_j)x} f_{\tau_j}(x), & \text{if } \tau_j = s_1, \\ e^{-u_1 k(\Upsilon,\tau_j)x} f_{\tau_j}(x), & \text{if } \tau_j \neq s_1, \tau_j \le \Upsilon. \end{cases}$$

(ii) *Given* $T_1 = t_1$, $S_1 = s_1 \notin M$ *and* $U_1 = u_1$, *the posterior parameters are* $M^\star = M \cup \{s_1\}$, *with*

$$f_{s_1}(x) = \mathcal{G}(x|1, \beta(s_1) + u_1 k(\Upsilon, s_1)),$$
$$f^\star_{\tau_j}(x) \propto e^{-u_1 k(\Upsilon,\tau_j)x} f_{\tau_j}(x) \qquad \text{if } \tau_j \le \Upsilon.$$

*Furthermore, given* $T_1 = t_1$ *and* $U_1 = u_1$, *the posterior distribution for the continuous part* $L^c_A(\cdot)$ *is* $L^c_A(s) \sim \Gamma_A\{\alpha(s), \beta^\star(s)\}$, *where* $\beta^\star(s) = \beta(s) + u_1 k(\Upsilon, s)$. *Thus, the posterior distribution of the normalized random measure is* $F^\star(t) = Z^\star(t)/\bar{Z}^\star$ *with* $Z^\star(t) = \int k(t,x)\,dL^\star_A(x)$.

Proposition 4 also holds for $\Upsilon = \infty$. However, for simulation purposes we need to truncate at $\Upsilon$. Given this result, posterior simulation becomes quite straightforward. For $n$ observations we have

$$p(t^n, s^n, u^n | L_A) = \prod_{i=1}^n \exp(-u_i \bar{Z}\,) k'(t_i, s_i)\, dL_A(s_i).$$

Given $L_A$, sampling from $p(s_i|t^n, u^n, L_A)$ and $p(u_i|t^n, s^n, L_A)$ is trivial, and given $(t^n, s^n, u^n)$, the conditional posterior of $L_A$ remains an additive process. We will need to implement a Gibbs sampler in the following way. Assuming that $M = \varnothing$, then initiate the algorithm by generating $u_i \sim \mathcal{G}(1,1)$, $\hat{s}_i \sim \mathcal{U}(0, t_i)$ for $i = 1, \ldots, n$, where $\mathcal{U}$ denotes the uniform distribution. For iterations $h = 1, \ldots, H$ do:

1. Generate $L^{(h)}_A$ from $p(L_A|t^n, s^{n^{(h-1)}}, u^{n^{(h-1)}})$ with the following specifications:

(a) The Lévy measure is given by

$$\nu_{\alpha,s}(dv) = dv \int_{(-\infty,s]} v^{-1} \exp\{-v\beta^\star(x)\}\alpha(dx),$$

where $\beta^\star(x) = \beta(x) + k(\Upsilon, x)\sum_{i=1}^n u_i^{(h-1)}$.

(b) The set of fixed jumps $M^{(h)} = \{s_1^{\star(h)}, \ldots, s_m^{\star(h)}\}$ is formed by all different $\{s_i^{(h-1)}\}$ with $r_j^{(h)}$, $j = 1, \ldots, m$, the number of $s_i^{(h-1)} = s_j^{\star(h)}$ for $i = 1, \ldots, n$.

(c) The distribution of the fixed jumps $L_A\{s_j^{\star(h)}\}$ is

$$f^{(h)}_{s^\star_j} = \mathcal{G}\left(r_j^{(h)}, \beta(s_j^{\star(h)}) + k(\Upsilon, s_j^{\star(h)}) \sum_{i=1}^n u_i^{(h-1)}\right).$$

2. Generate $s_i^{(h)}$ from $p(s_i|t^n, u^{n^{(h-1)}}, L_A^{(h)})$ for $i = 1, \ldots, n$ given by

$$p(s_j|t^n, u^{n^{(h-1)}}, L_A^{(h)}) \propto k'(t_i, s_i)\, dL_A^{(h)}(s_i) \mathbb{I}_{(-\infty, t_i)}(s_i).$$

3. Generate $u_i^{(h)}$ from $p(u_i|t^n, s^{n^{(h)}}, L_A^{(h)})$ for $i = 1, \dots, n$ given by

$$p(u_i|t^n, s^{n^{(h)}}, L_A^{(h)}) = \mathcal{G}\left(u_i \middle| 1, \int_0^{\Upsilon} k(\Upsilon, x)\, dL_A^{(h)}(x)\right).$$

Remark. In order to simulate from the continuous part of the posterior Lévy process $L_A^c(s)$, one option, which we employed, is to use the Ferguson and Klass (1972) algorithm. An alternative is the inverse Lévy method adopted by Wolpert and Ickstadt (1998). Both rely on approximations, making finite an infinite number of jumps. See Walker and Damien (2000) for the ideas.

3.1. *Numerical example.* Let us consider the case in which $L_A = \Gamma_A$ and $k(t, x) = \frac{1}{a}[1 - \exp\{-a(t - x)\}]\mathbb{I}_{[0,t]}(x)$ $(a \in \mathbb{R}^+)$, where $\mathbb{I}$ denotes the indicator function. This kernel has been motivated and used by Nieto-Barajas and Walker (2004). Thus $F$ is an MDP or, better, a *Dirichlet driven random probability d.f.* of the form

$$F(t) = \int_0^t [1 - \exp\{-a(t - x)\}]\, d\mathscr{D}_A(x)$$

and its corresponding random density is given by

$$f(t) = \int_0^t a \exp\{-a(t - x)\}\, d\mathscr{D}_A(x).$$

In this case one easily verifies that the arithmetic mean $\int_0^\infty t\, dF(t)$ can be written as $\int_0^\infty (x + \frac{1}{a})\, d\mathscr{D}_A(x)$; hence, its distribution is that of $\int_0^\infty x\, d\mathscr{D}_A(x)$ shifted by the factor $1/a$. The posterior density function of the mean, having observed $T^n = t^n$, is given by a slight modification of (8). Since the expression is difficult to deal with, we resort to our simulation algorithm, having set $\alpha(dx) = \mathbb{I}_{[0,5]}(x)\, dx$ and $a = 2$. We simulated $n = 100$ data points from a $\mathcal{G}(1, 1)$. Recall that the jumps of an IAP, when using the Ferguson and Klass algorithm, are simulated in a decreasing order according to their size. We truncated the number of jumps by calculating the relative error of a new jump and keeping only the jumps whose relative errors are greater than 0.0001. We ran the Gibbs sampling for 10,000 iterations with a burn-in of 1,000, keeping the last 9,000 simulations to obtain posterior summaries. Figure 1 presents the prior and posterior estimates of the normalized increasing process $F$. The prior estimate is placed away from the true d.f. and the posterior estimate follows very closely the true d.f., as expected.

In Figure 2 we can observe the prior and posterior distributions of the mean for $g(t) = t$. Due to the fact that the data were generated from a $\mathcal{G}(1, 1)$, the true value of the mean is 1. The prior distribution of the mean is situated away from the true value of the mean and has a large variance. The prior expected value of the mean is 2.74. On the other hand, the posterior distribution of the mean has a small variance and is concentrated around 1. The posterior expected value of the mean is 1.05.

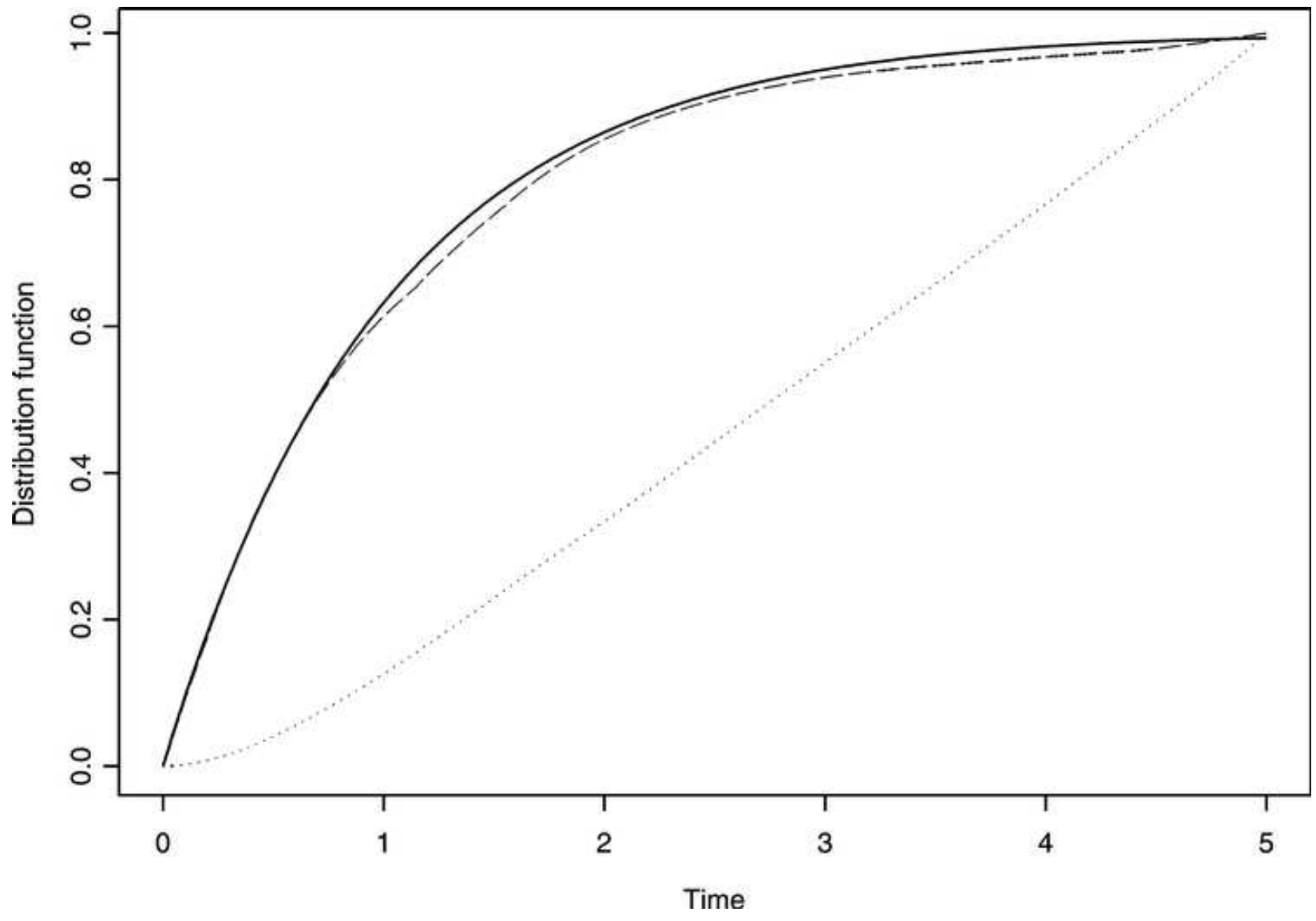

FIG. 1. *Prior and posterior estimates of the d.f. (———) True d.f., (········) prior estimate, (– – – – –) posterior estimate.*

## APPENDIX

**Details for the determination of conditions (I)–(III).** We have to show that, under (I)–(III), the sample paths of $F := \{F(t) = Z(t)/\bar{Z} : t \in \mathbb{R}\}$ are random probability d.f.s a.s. $[P]$.

Let us start with the denominator. We have to guarantee that $0 < \bar{Z} < +\infty$ a.s. $[P]$. Supposing $\int_{\mathbb{R}} \lim_{t\to+\infty} k(t,x)\, dL_A(x)$ is finite, we have $\bar{Z} = \int_{\mathbb{R}} \lim_{t\to+\infty} k(t,x)\, dL_A(x)$, a linear functional of a reparameterized IAP. Hence, Proposition 1 in RLP applies, leading one to state that $\bar{Z}$ is finite a.s. $[P]$ if and only if (II) holds.

Consider now the problem of the a.s. $[P]$ positiveness of $\bar{Z}$. Notice that, if (II) holds, we have

$$\exp\left\{-\int_{\mathbb{R}\times(0,+\infty)} [1 - \exp(-\lambda v \bar{k}(x))]\tilde{\nu}_\alpha(dx\, dv)\right\}$$
$$= E[e^{-\lambda \bar{Z}}]$$

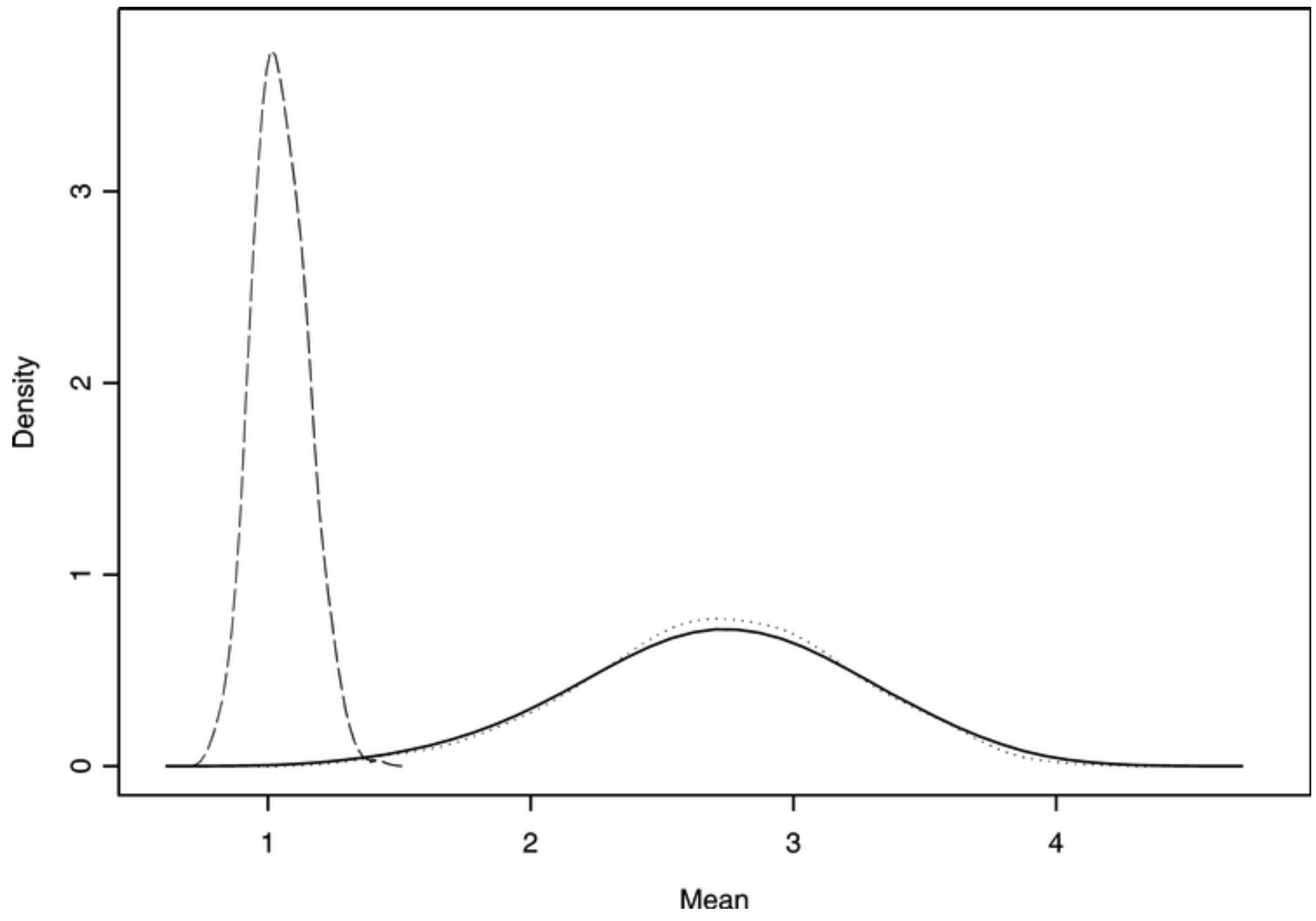

FIG. 2. *Prior and posterior distributions of the mean. (———) Exact prior distribution, (·········) simulated prior distribution, (– – – – –) simulated posterior distribution.*

$$= P\{\bar{Z} = 0\} + E[e^{-\lambda \bar{Z}} \mathbb{I}_{(0,+\infty)}(\bar{Z})],$$

where $\mathbb{I}$ denotes the indicator function. By the monotone convergence theorem, $P\{\bar{Z} = 0\} = \lim_{\lambda\to+\infty} \exp[\int_{-\mathbb{R}\times(0,+\infty)} [1 - \exp(-\lambda v \bar{k}(x))] \tilde{\nu}_\alpha(dx\, dv)]$. This entails that $P\{\bar{Z} = 0\} = 0$ if and only if $\lim_{\lambda\to+\infty} \int_{\mathbb{R}\times(0,+\infty)} [1 - \exp(-\lambda v \bar{k}(x))] \times \tilde{\nu}_\alpha(dx\, dv) = \infty$. Finally, we again apply monotone convergence so that $P\{\bar{Z} > 0\} = 1$ if and only if $\int_{\mathbb{R}\times(0,+\infty)} \tilde{\nu}_\alpha(dx\, dv) = +\infty$.

Turn to the numerator. In order to guarantee $t \mapsto Z(t)$ to be nondecreasing and right continuous it is enough to suppose $t \mapsto k(t, x)$ to be so for every $x \in \mathbb{R}$. Furthermore, if $\lim_{t\to-\infty} k(t, x) = 0$, we will have $\lim_{t\to-\infty} Z(t) = 0$.

PROOF OF PROPOSITION 3. (i) Let us start from the denominator. For every $s \in \mathbb{R}$,

$$E\left[\exp\left(is \int_{\mathbb{R}} \bar{k}(x)\, d\Gamma_A^\beta(x)\right)\right]$$

$$= \exp\left[-\int_{\mathbb{R}} \log\left(1 + is\frac{\bar{k}(x)}{\beta(x)}\right) dA(x)\right]$$

$$= \exp\left[-\int_{\mathbb{R}\times(0,+\infty)}\left(1-\exp\left(isv\frac{\bar{k}(x)}{\beta(x)}\right)\right)\frac{\exp(-v)}{v}\,dv\,dA(x)\right]$$

$$= E\left[\exp\left(is\int_{\mathbb{R}}\frac{\bar{k}(x)}{\beta(x)}\,d\Gamma_A\right)\right].$$

Applying the same arguments to the numerator, (i) follows.

(ii) The relation follows by application of the same arguments as in (i).

(iii) Follows immediately by application of Fubini's theorem. □

Proof of Theorem 1. In order to derive the posterior distribution, given $T^n = t^n$, we start by discretizing the MDP according to the procedure outlined in Section 2.2. The discretized random mean, at any fixed level $m$ of the tree of nested partitions, will be of the form $\int_{\mathbb{R}} g(t)\,dF_m(t) = \int_{\mathbb{R}}\sum_{j=0}^{k_m+1} g(a_{m,j})[k(t_{m,j+1},x)-k(t_{m,j},x)]\,d\mathscr{D}_A(x)$ and hence its d.f. can be written as

$$\begin{aligned}\mathbb{F}_m(\sigma) &= \frac{1}{2}\\ &\quad -\frac{1}{\pi}\int_0^{+\infty}\frac{1}{s}\,\mathrm{Im}\left\{\exp\left\{-\int_{\mathbb{R}}\log\left[1+is\left(\sum_{j=0}^{k+1} r_{m,j}(k(t_{m,j+1},x)\right.\right.\right.\right.\\ &\qquad\qquad \left.\left.-k(t_{m,j},x))-\sigma\right)\right]\\ &\qquad\qquad \left.\left.\times\,\alpha(dx)\right\}\right\}ds,\end{aligned}\tag{9}$$

where $r_{m,j} = g(a_{m,j})$, for $j = 0,\ldots,k_m+1$, and $\mathrm{Im}\,z$ stands for the imaginary part of $z \in \mathbb{C}$. Moreover, recall that

$$I^n_{a^+}h(\sigma) = \int_a^\sigma \frac{(\sigma-u)^{n-1}}{(n-1)!}h(u)\,du$$

is the Liouville–Weyl fractional integral, for $n \geq 1$, and $I^0_{a^+}$ represents the identity operator. By applying (5) to (9) together with some algebra, one obtains that its posterior density, given $T^n = t^n$, can be represented as

$$\rho_{t^n}(\sigma) = \begin{cases}\dfrac{(-1)^{q+1}}{\pi}I^{n-1}_{a^+}\,\mathrm{Im}\,\psi_m(\sigma), & \text{if } n = 2q,\\ \dfrac{(-1)^{q+1}}{\pi}I^{n-1}_{a^+}\,\mathrm{Re}\,\psi_m(\sigma), & \text{if } n = 2q+1,\end{cases}\tag{10}$$

with

$$\psi_m(\sigma) = \frac{1}{\mu^n(t^n)}$$

$$\times \int_0^{+\infty} s^{n-1} \int_{\mathbb{R}^n} \exp\Biggl(-\int_{\mathbb{R}} \log\Biggl[1 + is\Biggl(\sum_{j=0}^{k+1} r_j (k(t_{m,j+1}, x) - k(t_{m,j}, x)) - \sigma\Biggr)\Biggr] \times \alpha^*(dx)\Biggr) \times \prod_{p=1}^{n} [k(t_{m,i_p+1}, u_p) - k(t_{m,i_p}, u_p)] \times \alpha^n(du_1, \dots, du_n)\, ds,$$

where $\alpha^n$ is the $n$-fold product measure $\prod_{k=1}^n (\alpha + \sum_{i=1}^{k-1} \delta_{u_i})$ and $\alpha^*$ is given by $\alpha + \sum_{i=1}^n \delta_{u_i}$. In this case the expression for $\mu^n(t^n)$ is known, since it follows immediately by repeated application of Lemma 1 in Lo (1984),

$$\mu^n(t^n) = \Biggl(\prod_{i=1}^n (\alpha(\mathbb{R}) + i - 1)\Biggr)^{-1} \times \int_{\mathbb{R}^n} \prod_{p=1}^n [k(t_{m,i_p+1}, u_p) - k(t_{m,i_p}, u_p)] \alpha^n(du_1, \dots, du_n). \tag{11}$$

From (4) we know that (10) can be used as an approximate posterior density. Nevertheless, in this case we are able to obtain an explicit representation of the limiting posterior density.

Division of both numerator and denominator by $\prod_{p=1}^n [t_{m,i_p+1} - t_{m,i_p}]$ and application of the dominated convergence theorem yield the limiting posterior density which is given by (10) with

$$\psi(\sigma) = \frac{\prod_{i=1}^n (\alpha(\mathbb{R}) + i - 1)}{\int_{\mathbb{R}^n} \prod_{p=1}^n k'(t_p, u_p) \alpha^n(du_1, \dots, du_n)} \times \int_0^{+\infty} s^{n-1} \int_{\mathbb{R}^n} \exp\biggl\{-\int_{\mathbb{R}} \log[1 + is(h(x) - \sigma)] \alpha^*(dx)\biggr\} \times \prod_{p=1}^n k'(t_p, u_p) \alpha^n(du_1, \dots, du_n)\, ds. \tag{12}$$

Note that by Scheffé's theorem we have also convergence in total variation of $\mathbb{F}_{m,t^n}$ to $\mathbb{F}_{t^n}$. By application of Fubini's theorem it is possible to rewrite the posterior density function as

$$\rho_{t^n}(\sigma) = \int_{\mathbb{R}^n} \rho_{u^n}(\sigma) \frac{\prod_{p=1}^n k'(t_p, u_p) \alpha^n(du_1, \dots, du_n)}{\int_{\mathbb{R}^n} \prod_{p=1}^n k'(t_p, u_p) \alpha^n(du_1, \dots, du_n)}, \tag{13}$$

where, if $n = 2q$,

$$\rho_{u^n}(\sigma) = \frac{(-1)^{q+1} \prod_{i=1}^{n} (\alpha(\mathbb{R}) + i - 1)}{\pi} \times I_{a^+}^{n-1} \int_0^{+\infty} s^{n-1} \operatorname{Im}\biggl(\exp\biggl\{-\int_{\mathbb{R}} \log[1 + is(h(x) - \sigma)] \times \alpha^*(dx)\biggr\}\biggr) ds, \tag{14}$$

while, if $n = 2q + 1$, $\rho_{u^n}(\sigma)$ is obtained by simply substituting Im with Re in (14). Indeed, $\rho_{u^n}$ is a posterior density, given $U^n = u^n$, of a mean of a Dirichlet process, precisely of $\int_{\mathbb{R}} h(x)\, d\mathscr{D}_A(x)$. This can be seen by applying the procedure for derivation of posterior distributions of normalized RMI in Section 4 of RLP. Given the conjugacy of the Dirichlet process, we can replace (14) with the simpler expression

$$\rho_{u^n}(\sigma) = \frac{a}{\pi} \int_0^{+\infty} \operatorname{Re}\biggl(\exp\biggl\{-\int_{\mathbb{R}} \log[1 + is(h(x) - \sigma)]\alpha^*(dx)\biggr\}\biggr) ds. \tag{15}$$

Thus one has that a mean of an MDP, resulting from the combination of (13) and (15), is a mixture of a particular mean of a Dirichlet process, given the latent data $U^n$. By (11) it is easy to identify the mixing distribution as the distribution of $U^n$ conditionally on the real observations $T^n$. $\square$

Proof of Proposition 4. The idea of the proof is to express the likelihood function in a tractable way so we are able to apply standard Bayesian updating mechanisms. Let $T_1 = t_1$ be a single observation from $F$, and let $S_1 = s_1$ and $U_1 = u_1$ be auxiliary variables. Then the likelihood function is given by

$$\operatorname{lik}(L_A|t_1, s_1, u_1) = \exp\biggl\{-u_1 \int_0^{\infty} k(\Upsilon, x)\, dL_A(x)\biggr\} k'(t_1, s_1)\, dL_A(s_1).$$

Using product-integral properties [see, e.g., Gill and Johansen (1990)], the likelihood function can be rewritten as

$$\operatorname{lik}(L_A|t_1, s_1, u_1) = \biggl[\prod_{x \in [0,\infty)} \exp\{-u_1 k(\Upsilon, x)\, dL_A(x)\}\biggr] k'(t_1, s_1)\, dL_A(s_1).$$

Following Dykstra and Laud (1981), the prior process $L_A(\cdot)$ can be characterized by $dL_A^c(\nu) \sim \mathcal{G}(d\alpha(\nu), \beta(\nu))$ for the continuous part and $L\{\tau_j\} \sim f_{\tau_j}(x)$ for the prior fixed jumps. Based on the independence between increments in the prior process, the posterior conditional distribution for the continuous part and for the prior fixed jumps come straightforward. The only remaining point, to establish completely the posterior conditional distribution of $L_A(\cdot)$, is the distribution of the new fixed jump at $s_1$. For this, let

$$dL_A(s_1) = L_A[s_1, s_1 + \varepsilon),$$

and then

$$p(L_A[s_1, s_1+\varepsilon)|t_1, s_1, u_1) \propto L_A[s_1, s_1+\varepsilon)^{\alpha[s_1,s_1+\varepsilon)} e^{-\{u_1 k(\Upsilon,s_1)+\beta(s_1)\}L_A[s_1,s_1+\varepsilon)}.$$

Taking the limit as $\varepsilon \to 0$, we finally obtain that

$$p(L_A\{s_1\}|t_1, s_1, u_1) \propto e^{-\{u_1 k(\Upsilon,s_1)+\beta(s_1)\}L_A\{s_1\}},$$

as stated in the proposition. □

**Acknowledgments.** Special thanks are due to Antonio Lijoi for helpful suggestions and to Eugenio Regazzini for an early reading of the manuscript. The authors are grateful to an Associate Editor and two anonymous referees for their valuable comments that led to a substantial improvement in the presentation.

L. E. Nieto-Barajas
Departmento de Estadística
ITAM
Río Hondo 1
Col. Tizapán San Ángel
01000 México D.F.
México
e-mail: lnieto@itam.mx

I. Prünster
Dipartimento di Economia Politica
e Metodi Quantitativi
Università degli Studi di Pavia
Via San Felice 5
27100 Pavia
Italy
e-mail: igor.pruenster@unipv.it

S. G. Walker
Institute of Mathematics, Statistics
and Actuarial Science
University of Kent
Canterbury CT2 7NZ
United Kingdom
e-mail: S.G.Walker@kent.ac.uk